\documentclass[12pt,twoside]{amsart}

\date{}
\allowdisplaybreaks[4] \footskip=15pt
\renewcommand{\uppercasenonmath}[1]{}

\numberwithin{equation}{section} \theoremstyle{plain}
\newtheorem*{thm*}{Main Theorem}
\newtheorem{thm}{Theorem}[section]
\newtheorem{cor}[thm]{Corollary}
\newtheorem*{cor*}{Corollary}
\newtheorem{lem}[thm]{Lemma}
\newtheorem*{lem*}{Lemma}
\newtheorem{prop}[thm]{Proposition}
\newtheorem*{prop*}{Proposition}

\newtheorem*{que*}{Question}
\newtheorem{rem}[thm]{Remark}
\newtheorem*{rem*}{Remark}
\newtheorem{exa}[thm]{Example}
\newtheorem*{exa*}{Example}
\newtheorem{df}[thm]{Definition}
\newtheorem*{df*}{Definition}

\newtheorem*{conj*}{Conjecture}
\newtheorem*{ack*}{ACKNOWLEDGEMENTS}

\newcommand{\pf}{\noindent\begin {proof}}
\newcommand{\epf}{\end{proof}}

\begin{document}
\title{Strongly Goldie Dimension}
\author{L. Shen and J.L. Chen}

\address{Department of Mathematics, Southeast University, Nanjing, 210096,
P.R.China }
\address{Department of Mathematics, Southeast University, Nanjing, 210096,
P.R.China }
\date{May, 2005}
 \maketitle \baselineskip=20pt
\begin{abstract}
Let $R$ be an associative ring with identity. A unital right
$R$-module $M$ is called strongly  finite dimensional if
Sup$\{{\rm G.dim} (M/N)~|~\\ N\leq M\} < +\infty$. Properties of
strongly finite dimensional modules are explored. It is also
proved that: (1)If $R$ is left $F$-injective and strongly right
finite dimensional, then $R$ is left finite dimensional. (2) If
$R$ is right $F$-injective, then $R$ is right finite dimensional
if and only if $R$ is semilocal. Thus the Faith-Menal conjecture
is true if $R$ is strongly right finite dimensional. Some known
results are obtained as corollaries.
\end{abstract}

\bigskip

\bigskip
\section{Introduction}
Throughout this paper rings are  associative  with identity and
all modules are unital right $R$-modules. For a subset $X$ of a
ring $R$, the left annihilator of $X$ in $R$ is ${\bf l}(X)=\{r\in
R: rx=0$ for all $x\in X\}$. For any $a\in R$, we write ${\bf
l}(a)$ for ${\bf l}(\{a\})$. Right annihilators are defined
analogously. $N\leq_{ess} M$ means $N$ is an essential submodule of $M$.
Let $\mathcal A$ be a family of modules, we use $|\mathcal A|$ to denote
the cardinality of $\mathcal A$.\\
\indent An $R$-module $M$ has Goldie dimension $n$ (written G.dim
$M=n$) if there is an essential submodule $V\leq_{ess} M$ that is
a direct sum of $n$ uniform submodules. If, on the other hand, no
such integer $n$ exists, we write G.dim $M=+\infty$. We call an
$R$-module $M$ is finite dimensional if G.dim $M<+\infty$. A ring
$R$ is called right finite dimensional if it is finite dimensional
as a right $R$-module. Left dimensional rings can be defined
similarly. In this article, strongly Goldie dimension is
introduced. An $R$-module $M$ is with strongly Goldie dimension
$n$ (written SG.dim $M=n$) if Sup$\{{\rm G.dim} (M/N)~|~ N\leq
M\}=n$. Otherwise, SG.dim $M=+\infty$. $M$ is called strongly
finite dimensional if SG.dim $M<+\infty$. Properties of strongly
finite dimensional modules are explored. As applications, we show
that the Faith-Menal conjecture is true if $R_{R}$ is strongly
finite dimensional. The Faith-Menal conjecture was raised by Faith
and Menal in \cite{FM94}. It says that every strongly right Johns
ring is $QF$. Recall that a ring $R$ is called $QF$ if it is
one-sided noetherian and one-sided self-injective. A ring $R$ is
right Johns if $R$ is right noetherian and every right ideal is an
annihilator. Right Johns rings were characterized by Johns in
\cite{J77}, but he used a false result of Kurshan \cite[Theorem
3.3]{K70} to show that right Johns rings are right artinian. In
\cite{FM92}, Faith and Menal gave a counter example to show that
right Johns rings need not be right artinian. Later (see
\cite{FM94}) they defined strongly right Johns ring(the matrix
ring M$_{n\times n}$($R$) is right Johns for all $n \geq 1$) and
characterized such rings as right noetherian and left
$FP$-injective rings. But they didn't know whether a strongly
right Johns ring is $QF$.

\section{Characterizations of strongly finite dimensional modules}
\bigskip
\begin{df}\label{def 2.1}
{\rm An $R$-module $M$ is with strongly Goldie dimension $n$
(written SG.dim $M=n$) if Sup$\{{\rm G.dim} (M/N)~|~ N\leq M\}=n$.
Otherwise, SG.dim $M=+\infty$. $M$ is called strongly finite
dimensional if SG.dim $M < +\infty$. A ring $R$ is called strongly
right finite dimensional if it is strongly finite dimensional as a
right $R$-module. Strongly left finite dimensional rings are
defined similarly. }

\end{df}

\begin{exa}\label{exa 2.2}
If a ring $R$ is strongly right finite dimensional, then it is
right finite dimensional. But the converse is not true, even if
$R$ is a commutative  noetherian ring.
\end{exa}
\begin{proof}
For example, let $R$= $\mathbb{Z}$, then $R$ is a commutative
noetherian ring. But it is not strongly finite dimensional. For
every ideal of  $\mathbb{Z}$ is of the form $n \mathbb{Z}$ where
$n=p_{1}^{m_{1}}\cdots p_{k}^{m_{k}}$ is a product of powers of
prime numbers. Write $\mathbb{Z}_{l}= \mathbb{Z}/l\mathbb{Z}, l$
is a positive integer. Then it is well known that
$\mathbb{Z}_{n}\simeq
\mathbb{Z}_{p_{1}^{m_{1}}}\times\cdots\times\mathbb{Z}_{p_{k}^{m_{k}}}$,
which implies that G.dim $\mathbb{Z}_{n}\geq k$. Since $n$ is
arbitrary, SG.dim $\mathbb{Z}=+\infty$.
\end{proof}
\begin{prop}\label{prop 2.3}
If N is a quotient module of M, then SG.dim N $\leq$ SG.dim M.
\end{prop}
\begin{proof}
By definition.
\end{proof}
\begin{lem}{\rm (See \cite[Theorem 6.37]{L98})}\label{lem 2.4}
Let $0\rightarrow A\rightarrow B\rightarrow C\rightarrow 0$ be an
exact sequence of modules. Then  {\rm G.dim} B $\leq$ {\rm G.dim}
A {\rm+ G.dim} C.
\end{lem}
\begin{prop}\label{prop 2.5}
Let $0\rightarrow A\rightarrow B\rightarrow C\rightarrow 0$ be an
exact sequence of modules. Then ${\rm SG.dim}~ B\leq {\rm SG.dim}
~A+{\rm SG.dim}~ C$.
\end{prop}
\begin{proof}
For any $K\leq B$, by the above lemma, ${\rm G.dim} \frac B{K}\leq
{\rm G.dim} \frac {A+K}{K}+{\rm G.dim} \\ \frac
{B/K}{(A+K)/K}={\rm G.dim} \frac A{A\cap K}+{\rm G.dim} \frac
B{A+K}\leq {\rm SG.dim}~ A+{\rm SG.dim} \frac B{A}={\rm SG.dim}
~A+{\rm SG.dim} ~C$. Since $K$ is arbitrary, we have ${\rm
SG.dim}~ B\leq {\rm SG.dim} ~A+{\rm SG.dim} ~C$.
\end{proof}
\begin{cor}\label{cor 2.6}
{\rm SG.dim(}A{\rm +}B{\rm)}$~ \leq ~${\rm SG.dim} A{\rm +}{\rm
SG.dim} B
\end{cor}
\begin{proof}
Since $0\rightarrow A\rightarrow A+B\rightarrow B/A\cap
B\rightarrow 0$ is an exact sequence, by the above proposition,
SG.dim$(A+B)\leq {\rm SG.dim} ~A+{\rm SG.dim}\frac B{A\cap B}\leq
{\rm SG.dim}~ A+{\rm SG.dim}~ B$.
\end{proof}

\begin{cor}\label{cor 2.7}
${\rm SG.dim}(M_{1}\oplus M_{2}\oplus\cdots\oplus M_{n})$ = ${\rm
SG.dim} ~M_{1}+ {\rm SG.dim}~ M_{2}+\cdots+ {\rm SG.dim}~ M_{n}$.
\end{cor}
\begin{proof}
We only prove for $n=2$, the others are similarly. For any
quotient module $K_{i}$ of $M_{i}$, $i$=1,2. It is clear that
$K_{1}\oplus K_{2}$ is a quotient module of $M_{1}\oplus M_{2}$.
Then by definition, ${\rm SG.dim}~ M_{1}+ {\rm SG.dim}~ M_{2}\leq
{\rm SG.dim} (M_{1}\oplus M_{2})$.
 On the other side, ${\rm SG.dim}~ M_{1}+ {\rm SG.dim}~ M_{2}\geq
{\rm SG.dim} (M_{1}\oplus M_{2})$ by the above corollary.
\end{proof}
A nonzero module is called uniform if every nonzero submodule is
an essential submodule. A module is called uniserial if its
submodules are linearly ordered by inclusion. It is obvious that
every uniserial module is uniform. But the converse is not true.
For example, consider  $\mathbb{Z}$ as a $\mathbb{Z}$-module, then
$\mathbb{Z}$ is uniform but not uniserial. A ring $R$ is called
right serial if $R_{R}$ is a direct sum of uniserial modules.
\begin{cor}\label{cor 2.8}
If $R$ is right serial, then $R$ is strongly right finite
dimensional and SG.dim $R_{R}$=G.dim $R_{R}$.
\end{cor}
\begin{rem}\label{rem 2.9}
Since there are right artinian rings which are not right serial,
strongly finite dimensional rings may not be right serial.
\end{rem}
  \indent It is obvious that SG.dim $M\geq 1$ for
every nonzero module M. \\ If SG.dim $M= 1$, we have
\begin{thm}\label{thm 2.10}
The following are equivalent for an $R$-module $M$:\\
{\rm (1)} {\rm SG.dim} M {\rm= 1}.\\
{\rm (2)} $M$ is uniserial.\\
 {\rm (3)} Every nonzero quotient module of $M$ is uniform.\\
\end{thm}
\begin{proof}It is obvious that every quotient module of a uniserial module
is also uniserial. So $(2)\Rightarrow (3)\Rightarrow (1)$. Now we
assume (1), if $M$ is not uniserial, then there exists two
different submodules $A$ and $B$ of $M$, neither $A\subset B$ nor
$B\subset A$. Then it is clear that $\frac A{A\bigcap B}$ and
$\frac B{A\bigcap B}$ are nonzero submodules of $\frac M{A\bigcap
B}$. Hence there is an inclusion map: $\frac A{A\bigcap
B}$$\oplus$$\frac B{A\bigcap B}\hookrightarrow$$\frac M{A\bigcap
B}$, which implies SG.dim $M\geq 2$, a contradiction.
\end{proof}
\begin{exa}\label{exa 2.11}
A strongly finite dimensional module M may not be noetherian, even
if {\rm SG.dim} M {\rm= 1}.
\end{exa}
\begin{proof}
For example, let $p$ be a positive prime number. Then
$M=\{a/p^{n}\in \mathbb{Q} ~| ~a\in \mathbb{Z}$ and
$n\in\mathbb{N}\}$ is an additive subgroup of $\mathbb{Q}$ with
subgroup $\mathbb{Z}$. Denote the factor group $M/\mathbb{Z}$ by
$\mathbb{Z}_{p^{\infty}}$. It is clear that every subgroup of
$\mathbb{Z}_{p^{\infty}}$ is cyclic and spanned by $1/p^{n}$ for
some $n$, which implies that $\mathbb{Z}_{p^{\infty}}$ is a
uniserial but not noetherian $\mathbb{Z}$-module.
\end{proof}
\begin{exa}\label{exa 2.12}
{\rm A strongly finite dimensional module M may not be artinian,
even if {\rm SG.dim} M {\rm= 1}. For example (see \cite[Example
2.12]{P01}), let $\mathbb{Z}[i]$ be the ring of Gaussian's
integers, $V=\mathbb{Z}[i]_{(2-i)}$ be its localization at the
prime ideal generated by $2-i$, $\sigma$ be the complex
conjugation and $R$ be the skew power series ring
$R=\{\alpha_{0}+x\alpha_{1}+x^{2}\alpha_{2}+\cdots
~|~\alpha_{0}\in V, \alpha_{j}\in\mathbb{Q}[i] ~for j\geq 1\}$,
where the multiplication is given by the rule $\alpha
x=x\sigma(\alpha)$. Then $R$ is uniserial but not artinian.}
\end{exa}

Let $M$ be an $R$-module and $\{K_{\lambda}\}_{\wedge}$ be a
family of proper submodules of $M$. $\{K_{\lambda}\}_{\wedge}$ is
called coindependent (see \cite{L99}) if for every
$\lambda\in\wedge$ and finite subset $I\subseteq \wedge \setminus
\{\lambda\}$ $K_{\lambda}+\cap_{i\in I}K_{i}=M$ (if $I$ is the
empty set, then set $\cap_{i\in I}K_{i}=M$).
\\  For SG.dim $M=n\geq 2$, we have
\begin{thm}\label{thm 2.13}
The following are equivalent for an R-module $M$:\\
{\rm (1)} {\rm SG.dim}$M$ {\rm = n}.\\
{\rm (2)} {\rm Sup} $\{|\mathcal A|~|~\mathcal A$ is a
coindependent family of N, where N is any submodule of a quotient
module of $M \}=n$
\end{thm}
\begin{proof}
(1)$\Longrightarrow$ (2). Assume there exists a nonzero submodule
$K$ of a quotient module  of $M$ and a coindependent family
$\{K_{i}, 1\leq i\leq m\}$ of $K$ such that $m\geq n+1$. Then we
have a canonical inclusion $f: \frac
K{\cap_{\i=1}^{m}K_{i}}\hookrightarrow \frac K K_{1}\oplus
\cdots\oplus\frac K K_{m}$ with
$f(k+\cap_{\i=1}^{m}K_{i})=(k+K_{1},\cdots, k+K_{m}), \forall k\in
K$. Now we set $N_{i}=\cap_{j\neq i}K_{j}, 1\leq i\leq m$. Then by
the definition of coindependent family, $K_{i}+N_{i}=K, 1\leq
i\leq m$. Thus for each $i$, $f(N_{i})$ is not zero and
$f(N_{i})=(0,\cdots,x+K_{i}, \cdots,0),~ x\in N_{i}$. Hence
$f(N_{1})+\cdots +f(N_{m})$ is a direct sum. Since $f$ is monic,
G.dim$(\frac K{\cap_{\i=1}^{m}K_{i}})\geq m$, which shows that
G.dim $(\frac M{\cap_{\i=1}^{m}K_{i}})\geq m>n$, a
contradiction.\\
(2)$\Longrightarrow$ (1). If SG.dim $M > n$, then there exists a
nonzero quotient module $N$ of $M$ with a family $\{K_{i}, 1\leq
i\leq m>n\}$ of independent submodules of $N$. Now take
$N_{i}=\sum_{k\neq i}K_{i}, 1\leq i\leq m$. It is clear that
$\{N_{i}, 1\leq i\leq m\}$  is a coindependent family of $\sum
K_{i}$, which is a contradiction. So SG.dim $M \leq n$. By
hypothesis, there exists a nonzero module $A$, which is a
submodule of a quotient module $B$ of $M$. And $A$ has a
coindependent family $\{A_{i}, 1\leq i\leq n\}$ of proper
submodules of $B$. Then by the construction in the proof of
(1)$\Longrightarrow$ (2), G.dim $\frac A{\cap_{\i=1}^{n}A_{i}}\geq
n$. So SG.dim $M\geq$ G.dim $\frac
B{\cap_{\i=1}^{n}A_{i}}\geq$G.dim $\frac
A{\cap_{\i=1}^{n}A_{i}}\geq n$. Hence SG.dim $M=n$.
\end{proof}

A module $M$ is said to be of finite length (see \cite{AF92}) if
there exists a finite chain of submodules of $M$:
$M=M_{0}>M_{1}\cdots > M_{n}=0$ where $M_{i-1}/M_{i}$ is simple
$(i=1,2,\ldots,n)$. Write $c(M)=n$ for the composition length of
$M$.
\begin{lem}{\rm [The Schreier Refinement Theorem]}\label{lem 2.14}\\
If M is a module of finite length and if $M=N_{0}>N_{1}\cdots >
N_{p}=0$ is a chain of submodules of M, then there is a
composition series for M whose terms include $N_{0}, N_{1}, \cdots
, N_{p}=0$
\end{lem}
\begin{lem} \cite[Corollary 6.7 (2)]{L98}\label{lem 2.15}
If $M$ is a module of finite length, then $G.dim~ M\leq c(M)$.
$G.dim ~M= c(M)<+\infty$ if and only if $M$ is semisimple.
\end{lem}
\begin{prop}\label{prop 2.16}
If $M$ is a module of finite length, then $SG.dim~M\leq c(M)$.
$SG.dim~M= c(M)<+\infty$ if and only if $M$ is semisimple.
\end{prop}
\begin{proof}
For any quotient module $N$ of $M$, $c(N)\leq c(M)$ by Lemma 2.14.
Thus by Lemma 2.15, SG.dim $M\leq c(M)$. If $M$ is a semisimple
module of finite length, then it is easy to get that SG.dim $M$=
$c(M)<+\infty$. On the converse, if SG.dim $M= c(M)<+\infty$, then
there exists a quotient module $N_{0}$ of $M$ such that G.dim
$N_{0}$=SG.dim $M$=$c(M)$. But by Lemma 2.14, $c(N_{0})<c(M)$ if
$N_{0}$ is not isomorphic to $M$. Thus $N_{0}$ must be isomorphic
to $M$ and G.dim $N_{0}$=$c(N_{0})$. Hence $M\simeq N_{0}$ is
semisimple by Lemma 2.15.
\end{proof}
A ring $R$ is called right semiartinian if every right $R$-module
has an essential right socle.
\begin{thm}\label {thm 2.17}
A ring R is right artinian if and only if R is right
semiartinian and strongly right finite dimensional.
\end{thm}
\begin{proof}
If $R$ is right artinian, then it is right semiartinian. Since any
right artinian ring is right noetherian,  $R_{R}$ is of finite
length by \cite[Proposition 11.1]{AF92}. Thus by the above
proposition, $R$ is strongly finite dimensional. On the other
hand, if $R$ is right semiartinian and strongly right finite
dimensional, then every cyclic right $R$-module is finitely
cogenerated. So $R$ is right artinian by V$\acute{a}$mous Lemma
(see \cite[Lemma 1.52]{NY03}).
\end{proof}
\begin{rem}\label {rem 2.18}
From Example 2.12, we see that right semiartinian condition is
necessary in the above theorem.
\end{rem}
\begin{exa}\label {exa 2.19}
{\rm There exists strongly right finite dimensional rings which
are not strongly left finite dimensional. Since every one sided
artinian ring is left and right semiartinian, by Theorem 2.15, we
only need to find right artinian rings which are not left
artinian.}
\end{exa}
\section{Applications}
Recall that a ring $R$ is called right $F$-injective if, for every
$R$-homomorphism from a finitely generated right ideal to $R_{R}$
can be extended from $R_{R}$ to $R_{R}$.  $R$ is called right
$FP$-injective if, for any free right $R$-module $F$ and any
finitely generated $R$-submodule $N$ of $F$, every
$R$-homomorphism $f: N\rightarrow R$ can be extended to an
$R$-homomorphism $g: F\rightarrow R$. It is obvious that right
$FP$-injective rings are right $F$-injective. But it is still
unknown whether a right $F$-injective ring is right
$FP$-injective. Left $F$-injective and left $FP$-injective rings
can be defined similarly.
\begin{lem}\label{lem 3.1}
{\rm (see \cite[Lemma 1.37]{NY03})}A ring R
is left F-injective if and only if it satisfies the following
conditions:\\
{\rm (1)} ${\bf r}(T\cap T^{'})={\bf r}(T)+{\bf r}(T^{'})$ for all
finitely generated left ideals T and $T^{'}$ of R.\\
{\rm (2)} ${\bf rl}(a)=aR$ for all $a\in R$.
\end{lem}
\begin{thm}\label{thm 3.2}
If R is left F-injective and strongly right finite dimensional,
then R is left finite dimensional.
\end{thm}
\begin{proof}
If $R$ is not left finite dimensional, then there exist $a_{i}\in
R, i=1,2,3,\cdots,$ such that $\sum_{i=1}^{\infty} Ra_{i}$ is a
direct sum. Since $R$ is left $F$-injective, by the above lemma,
for any integers $m$ and any finite subset $I\subset
\mathbb{N}\setminus m$ , $R={\bf r}(0)={\bf r}(Ra_{m}\cap
\sum_{i\in I}Ra_{i})$=${\bf r}(a_{m})+{\bf r}(\sum_{i\in
I}Ra_{i})={\bf r}(a_{m})+\cap_{i\in I}{\bf r}(a_{i})$.  It is also
clear ${\bf r}(i)\neq{\bf r}(j)$, $i\neq j$. Thus $\{{\bf
r}(a_{i}), i=1,2,3,\cdots\}$ is a coindependent family of $R_{R}$.
 So by Theorem 2.13, SG.dim$(R_{R})=+\infty$, a
contradiction.
\end{proof}
\begin{thm}\label{thm 3.3}
The Faith-Menal conjecture is true if R is strongly right finite
dimensional.
\end{thm}
\begin{proof}
Since $R$ is strongly right Johns, $R$ is left $FP$-injective. So
$R$ is left finite dimensional by the above theorem. Thus $R$ is
$QF$ by \cite[Corollary 1.3]{FM94}.
\end{proof}
\begin{thm}\label{thm 3.4}
If R is right F-injective, then R is semilocal if and only if R is
right finite dimensional.
\end{thm}
\begin{proof}
If $R$ is semilocal, then by \cite[Theorem 1.3, Corollary
3.2]{L99}, $_{R}R$ does not contain an infinite coindependent
family of submodules. Then from the proof in Theorem 3.2, we have
that $R$ is right dimensional. On the converse, since $R$ is right
$F$-injective, $R$ is right $C2$ (every right ideal is essential
in a direct summand of $R_{R}$) by \cite[Proposition 5.10]{NY03}.
Thus $R$ is semilocal by \cite[Corollary C.3]{NY03}.
\end{proof}
\begin{cor}\label{cor 3.5}
If $R$ is right FP-injective or right self-injective, then R is
semilocal if and only if R is right finite dimensional.
\end{cor}
\begin{cor}\label{cor 3.6}
{\rm (see \cite[Theorem 5.56]{NY03})}The following are equivalent
for a ring R: \\
{\rm (1)} R is semilocal, right FP-injective, and right Kasch.\\
{\rm (2)} R is right finite dimensional, right FP-injective, and
right Kasch.
\end{cor}
A ring $R$ is called right $PF$ if $R$ is a semilocal and right
self-injective ring with an essential right socle.
\begin{cor}\label{cor 3.7}
R is right PF if and only R is a right finite dimensional and
right self-injective ring with an essential right socle.
\end{cor}
\begin{center}
{ACKNOWLEDGEMENTS}
 \end{center}
\indent\indent The research is supported by the National Natural
Science Foundation of China (No.10171011) and the Teaching and
Research Award Program for Outstanding Young Teachers in Higher
Education Institutes of MOE, P.R.C.


\begin{thebibliography}{99}
\bibitem{AF92} F. W. Anderson and K. R. Fuller, Rings and Categories of Modules,
 Second Edition. Graduate Texts Vol. {\bf 13} Springer-Verlag. Berlin-Heidelberg-New York. 1992.

\bibitem{FM92}Carl Faith and Pere Menal, A counter example to a conjecture of Johns,
{\it Proc. A.M.S}, {\bf 116} {1992}, 21-26.
\bibitem{FM94}Carl Faith and Pere Menal, The Structure of Johns
Rings, {\it Proc. A.M.S}, {\bf 120} {1994}, 1071-1081.

\bibitem{J77} Baxter Johns, Annihilator Conditions in Noetherian
Rings, {\it J. Algetra}, {\bf 49} (1977), 222-224.
 \bibitem{K70} R.P. Kurshan, Rings whose cyclic modules have finitely genrated socle,
 {\it J. Algetra}, {\bf 15} (1970), 376-386.
\bibitem{L98} T. Y. Lam, Lectures on Modules and Rings.  Graduate Texts in Mathematics
vol{\bf 189}. Springer-Verlag. New York. 1998.

\bibitem{L99} C. Lomp,  On semilocal modules and
rings, {\it Comm. Algebra}, {\bf 27(4)} (1999), 1921-1935.


\bibitem{NY03} W. K. Nicholson,  M. F. Yousif,  {\it Quasi-Frobenius
Rings}, Cambridge Tracts in Mathematics, 158, Cambridge University
Press, Cambridge, 2003.
\bibitem{P01} G. Puninski, Serial Rings, Kluwer Academic
Publishers, Dordrecht/ Boston/ London, 2001.







\end{thebibliography}
\end{document}